\documentclass[11pt]{article}
\usepackage{graphicx}
\usepackage{eucal}
\usepackage{amssymb}
\newtheorem{theorem}{Theorem}
\newtheorem{proposition}[theorem]{Proposition}

\newtheorem{definition}[theorem]{Definition}
\newtheorem{example}[theorem]{Example}
\newtheorem{remark}[theorem]{Remark}
\newtheorem{lemma}[theorem]{Lemma}

\def\bM{{\bf M}}
\def\bB{{\bf B}}

\def\expo{{\rm expo}}
\def\End{{\rm End}}

\def\ra{\longrightarrow}
\def\da{\downarrow}

\def\End{{\rm End}}

\def\rank{{\rm rank}}

\def\ker{{\rm ker}}
\def\im{{\rm im}}
\def\whM{\widehat{M}}

\def\whJ{\widehat{J}}
\def\whT{\widehat{T}}

\def\da{\downarrow}

\def\cF{\mathcal{F}}

\def\Gr{{\rm Gr}}
\def\C{{\mathfrak{C}}}

\title{The commutator algebra of a nilpotent matrix and an application to the theory of commutative Artinian algebras}

\author{Tadahito Harima   \\ 
Department of Mathematics \\
Hokkaido University of Education  (Kushiro Campus)\\
Kushiro  085-8580, Japan \\
Email: harima@kus.hokkyodai.ac.jp 
\and 
Junzo Watanabe  \\ 
Department of Mathematics, Tokai University \\ 
Hiratsuka 259-1292, Japan \\
Email: junzowat@keyaki.cc.u-tokai.ac.jp
}

\begin{document}
\maketitle


\begin{abstract}
We show a number of properties of the commutator algebra
of a nilpotent matrix over a field.  In particular we determine the
simple modules of the commutator algebra.
Then the results are applied to prove that certain
Artinian complete intersections have the strong
Lefscehtz property.
\end{abstract}

\tableofcontents

\section{Introduction}
Suppose that $A =\bigoplus ^c  _{i=0}A_i$ is an Artinian graded algebra over a field $K$ 
and $z$ is a homogenous element of $A$.   
We denote by $\times z$ the linear  map $A \ra A$ defined  by $\times z(a)=az$.  
We say that  $A$ has the weak Lefschetz property if there is a linear  form 
$g \in A$ such that the rank of $\times g$ is maximum that can be expected only from the 
Hilbert function of $A$.  In the same sense the strong Lefschetz property means that the rank of 
$\times g^k$ is maximum possible for every $k \geq 0$.  
(For details see Definition~\ref{def_cosperner}).  We use the definition of 
``strong Lefschetz property'' in a restricted sense.  
For the reason see \cite[Remark~4]{tHjW03}. 
Consider the exact sequence 
$$0 \ra A/(0:z) \ra A \ra A/(z) \ra 0$$
where the first map is 
induced  by the multiplication map $\times z$.  
It seems to be interesting to ask under what conditions the strong/weak  
Lefschetz property of $A$ can be deduced from that of $A/(0:z)$ and $A/(z)$. 
Let $g$ be a general linear form (one independent of $z$).  For the definition of 
``general linear form,'' see Definition~\ref{def_cosperner}.    
To determine the rank of the linear  map $\times g : A \ra A$ from the knowledge of 
the two homomorphisms 
$\times g : A/(z) \ra A/(z)$ and $\times g: A/(0:z) \ra A/(0:z)$, one has to consider how a preimage 
of $A/(z)$ in $A$ is mapped to $(z)$ by the map $\times z$. 
Thus we are led to consider the exact sequence 
$0 \ra (0:z) \ra A \ra (z) \ra 0$ as well.  
Now to determine the rank of the linear map 
$$\times z : A/(0:z) \ra A/(0:z)$$ we consider 
the exact sequences
$$0 \ra A/((0:z):z) \ra A/(0:z) \ra A/(0:z)+(z) \ra 0$$ 
and
$$0 \ra (0:z):z \ra A/(0:z) \ra (z)+(0:z)/(0:z) \ra 0.$$  
These sequences  may be used repeatedly, so as to complete the calculation of the rank of $\times z$. 

In our paper \cite{tHjW06}, we consider a family ${\mathcal F}$
of $K$-algebras  
 such that if $A \in \cF$ then there is an exact sequence 
 $0 \ra A' \stackrel{\times z}{\ra} A \ra B \ra 0$, with some linear form $z$, 
 where 
 $A' \in \cF$ and $B \in \cF$,  and the rank of $\times g : A \ra A$ for a 
 general linear form  can be determined from the knowledge of $A'$ and $B$.
The purpose of this paper is to provide basic tools and framework for considering 
such a family of algebras, and to prove the strong/weak Lefschetz property for all 
members of the family altogether.  

One of the basic tools is the commutator algebra of a nilpotent matrix.  
A famous lemma due to I. Schur  states that the commutator 
algebra of an irreducible matric set 
is a division ring.   
In some sense we are interested in the other extremal case,  
namely the commutator algebra of a single nilpotent matrix. 
In our previous paper~\cite{tHjW03} the commutator algebra of some commuting family of 
matrices in a full matrix algebra played an important role in studying  commutative Artinian algebras. 
In this paper we treat the idea more generally. 
For a nilpotent matrix $J$ we denote by $\C(J)$ the commutator algebra of $J$.
A description of  $\C(J)$ can be found in \cite{fG59},  \cite{hTaA32},  and \cite{rB03}. 
The algebra $\C(J)$ as a set of matrices is easily determined as it is the set of solutions of a system of linear equations.   
The way  it is put as a set of matrices, however,  is not quite adequate for our purposes.
We  rearrange the order of the basis elements,  so that the matrices of the commutator algebra 
are ``upper triangular'' as much as possible.   The resulting set of matrices is denoted by 
$\C(\whJ)$.  The use of the transformation  
$J \mapsto \whJ$  was suggested by  \cite{hW46}  and it worked  quite effectively in \cite{tHjW03}.
We show the relation between 
$\C(\whJ)$ and $\C(\whJ ')$, where $\whJ'$ is the submatrix of $\whJ$ 
corresponding to the restricted map 
$J: \im\; J \ra \im\; J$.  This is important in the 
inductive argument and it reveals the 
structure of the algebra $\C(J)$. 
Once this is done, we may determine the Jacobson radical of the commutator algebra;  
we also determine in several ways the simple modules of $\C(J)$ 
(Proposition~\ref{central_module}). 
In our papers \cite{tHjW06} and \cite{tHjW07} we study further these central simple 
modules for Gorenstein algebras $A$ with a linear form fixed.    

Furthermore we show how the rank of 
$\whM +  \whJ$ can be  computed (Proposition~\ref{rank_of_big_matrix}) for a certain 
element $\whM \in \C(\whJ)$.  
In our application to commutative Artinian  $K$-algebras  
we are going to choose two linear elements,  $l, z$  of $A$,  and show that, with certain conditions 
imposed on $A$ and $z \in A$, 
the rank of $\times (l + \lambda z)$, for most $\lambda \in K$,  
reaches the CoSperner number of $A$, so the pair $(A, l + \lambda z)$ is weak Lefschetz.  

The main results of this paper are Proposition~\ref{rank_of_deformation_2} 
and Theorem~\ref{rank_of_general_element} of Section~4. 
To explain the meaning of  Theorem~\ref{rank_of_general_element}, 
let $y, z \in A$ be linearly independent linear forms of 
an Artinian algebra $A$.  The theorem gives a lower bound for the rank 
$$\times (y +  \lambda z) \in \End(A)$$ 
in terms of $\dim A/(z)$ and the ranks of the diagonal blocks of the linear maps   
$$\times y ^i : A/(z) \ra A/(z) \in \End(A/(z)), \ i=1,2,\cdots.  $$
Since the obvious upper bound for the rank $\times (y + \lambda z)$ is the 
CoSperner number of $A$,  these considerations give us a sufficient condition for $A$ to have the 
weak Lefschetz property (Theorem~\ref{rank_of_general_element} $(ii)$).    
This is a direct consequence of Propositions~\ref{rank_of_big_matrix} 
and \ref{rank_of_deformation_2}, since $\times y \in \C(\times z)$. 
To prove Proposition~\ref{rank_of_deformation_2},  
we consider the form ring 
$$\Gr_{(z)}(A)=A/(z)\oplus (z)/(z^2) \oplus (z^2)/(z^3) \oplus \cdots \oplus (z^{p-1})/(z^p)$$
for a $K$-algebra $A$ with respect to the principal ideal $(z)$.  
It is well known that this is endowed with algebra structure; 
in fact it is isomorphic to an algebra  $B[y]/I$ where  $B=A/(z)$ and    
$I$ is generated by ``monomials in $y$,''  which are elements of the form  $\overline{a_i}y^{m_i}$.  
It turns out that the inequality of the proposition is essentially equivalent to 
$$\rank (\times g_1) \leq \rank (\times g_2)$$ 
where
$g_1$ and $g_2$ are general linear forms  of $\Gr_{(z)}(A)$ and  $A$ and  respectively. 
This simplifies the computation of the rank of
a general linear form  as it is often the case that  even 
the rank of $\times g_1$ is sufficiently large.     

Throughout this paper $K$ denotes a field of any characteristic unless otherwise specified.  
When we discuss the strong Lefschetz property  we assume the characteristic of $K$ is zero. 
Except for this all results are valid for any characteristic. 
In  Section~2.2  and an early part of Section~2.3 we do use the exponential of a nilpotent matrix.  
However, this is {\em not} essential.  We use the exponential, because it facilitates  
the description of the matrices of the commutator algebra $\C(J)$. This is easy to see, so 
there will be no confusion.   In \cite{tHjW07}, we consider Artinian algebras with non-standard grading, but 
in this paper all algebras are assumed to have  the standard grading so the generators  of algebras have degree one. 

The authors would like to express their thanks to the referee for invaluable suggestions for this paper.  
\section{The commutator algebra of a nilpotent matrix }  
\subsection{The Jordan canonical form and Young diagrams}
Let  $K$  be  a field.  
Let    $\bM(n)$    denote  the set of $n \times n$  matrices with entries in  $K$.
We are concerned with nilpotent matrices.  
Since all eigenvalues of a nilpotent matrix are 0,  its  Jordan canonical form  is 
expressed by telling  how it decomposes into Jordan cells.   
Thus we are led to use a Young diagram to indicate the Jordan canonical form of 
a nilpotent matrix.    When we say that $(n_1, n_2, \cdots, n_r)$   
is a partition of an integer $n$, it means that all $n_i$ are positive integers such that 
$n=n_1+n_2+ \cdots + n_r$.  When the terms in a partition of $n$ 
are arranged in decreasing order we may associate to it a Young diagram of 
size $n$ in the well known manner. 
The notation  $T=T(n_1, n_2, \cdots, n_r)$ denotes the Young diagram 
with rows of $n_i$ boxes, $i=1,2,\cdots, r$, where it is tacitly assumed that 
$n_1 \geq n_2 \geq \cdots \geq n_r >0$.  
The same Young diagram is also denoted by 
$T=\whT(\nu _1, \cdots, \nu _p)$, by which it is meant that 
the integer $\nu _i$ is the number of boxes in the $i$th column of $T$.

Let $T$ be a Young diagram of size $n$. 
Suppose we number the boxes of $T$ with numbers $1,2, \ldots, n$.  Then we may define a matrix  
$M=(a_{ij}) \in \bM(n)$ as follows:
\begin{center}
\begin{equation}  \label{jordan_matrix}
a_{ij}=\left\{ \begin{array}{l}   \label{def_of_nilpotent_mat}
1 \ \  \mbox{if $j$ is next to and on the right of $i$,}\\
0 \ \  \mbox{otherwise.}
\end{array}
\right. 
\end{equation}
\end{center}
It is easy to see that the matrix $M$ is nilpotent and the matrices given by 
different numberings of the same $T$ are different only by a permutation of 
rows and columns.  
As is well known there is a bijection between the set of Young diagrams of size $n$ and
the conjugacy classes of matrices  with a single eigenvalue in $\bM(n)$. 

In this note we use only two numberings of a Young diagram. 
One is the horizontal numbering starting with the first row and ending with the last row,   
and the other vertical.   

Here are  examples of such numberings for  the Young diagram $T=T(5,3,1)$. 
\newcommand{\ssvline}[1]{\multicolumn{1}{|c}{#1}}
\newcommand{\svline}[1]{\multicolumn{1}{|c|}{#1}}

\begin{center}
\begin{equation}
\begin{array}{|ccccc}                                                                 \hline
1   &  \svline{2}   & \svline{3}  & \svline{4}  & \svline{5}     \\   \hline
6   &  \svline{7}   & \svline{8}  &                   &                     \\    \cline{1-3}
9   & \ssvline{}    &                   &                   &                     \\    \cline{1-1}
\end{array} , 
\hspace{8ex}
\begin{array}{|ccccc}                                                                 \hline
1   &  \svline{4}   & \svline{6}  & \svline{8}  & \svline{9}     \\   \hline
2   &  \svline{5}   & \svline{7}  &                   &                     \\    \cline{1-3}
3   & \ssvline{}    &                   &                   &                     \\    \cline{1-1}
\end{array} 
\end{equation}
\end{center}

When $T$ is numbered vertically we call the matrix defined by (\ref{jordan_matrix}) the 
{\bf Jordan second canonical form}.
To write the matrix explicitly, let $n=\nu_1 + \cdots + \nu _p$ be the dual partition of $n=n_1+ \cdots + n_r$. 
In other words $\nu _i$ is the number of boxes of the $i$th column of $T(n_1, \cdots, n_r)$.
Let $I_i$ be the  $\nu_i \times \nu_{i+1}$  matrix 
\begin{equation}
I_{i}=\begin{array}{cl}                                               \\   \cline{1-1}
\svline{E}               &             \} \nu_{i+1}                  \\   \cline{1-1}
\svline{O}               &            \} \nu_{i}-\nu_{i+1}      \\  \cline{1-1}
  \underbrace{}_{\nu_{i+1}}         
\end{array}                    
\end{equation}
where $E$ is the  $\nu_{i+1} \times \nu_{i+1}$  identity and $O$ the $(\nu_{i}-\nu_{i+1}) \times \nu_{i+1}$  zero matrix.  
Then the Jordan second canonical form for $T=T(n_1, n_2, \cdots, n_r)$  is the matrix
\begin{equation}  \label{jordan_2nd_canonical_form}
\begin{array}{ccccccc}                                                                                                                                                              \\   \cline{1-6}
\svline{O}                   &  \svline{I_{1}}             &  \svline{O}              & \svline{\cdots}  &\svline{O}& \svline{O}              & \} \nu_1      \\   \cline{1-6}
\svline{O}                    &  \svline{O}                   &  \svline{I_{2}}        & \svline{\cdots}  &\svline{O}& \svline{O}                & \} \nu_2    \\  \cline{1-6}
\svline{O}                    &  \svline{O}                   &  \svline{O}               & \svline{\cdots}  &\svline{O}& \svline{O}                & \} \nu_3   \\  \cline{1-6}
\svline{\vdots}            &  \svline{\vdots}           &  \svline{\vdots}        & \svline{\ddots}  &                & \svline{\vdots}       &  \vdots     \\    \cline{1-6}
\svline{O}                    & \svline{O}                     & \svline{O}                 & \svline{\cdots}  &\svline{O}& \svline{I_{p-1}}      & \} \nu_{p-1}   \\    \cline{1-6}
\svline{O}                    & \svline{O}                     & \svline{O}                  & \svline{\cdots}  &\svline{O} & \svline{O}               & \} \nu_p   \\    \cline{1-6}
\underbrace{}_{\nu_1} &\underbrace{}_{\nu_2}  &\underbrace{}_{\nu_3} & \ldots               &                  & \underbrace{}_{\nu_p}
\end{array}                    
\end{equation}

When $T$ is numbered horizontally the matrix is usually called the Jordan canonical form.  
If it is necessary to distinguish from the second, we call it the Jordan first canonical form. 
In our previous paper \cite{tHjW03} we used the term ``Jordan second canonical form'' 
in a slightly different sense from the use here.

For our convenience we state  some notational conventions.
\begin{itemize}
\item
$T=T(n_1, n_2, \cdots, n_r$) indicates that $T$ is a Young diagram with $r$ rows, where $n_i$ is the number of 
boxes of the $i$th row.   It is assumed that 
$n_1 \geq n_2 \geq \cdots \geq  n_r > 0$.
$T=\whT(\nu _1, \cdots, \nu _p)$ denotes the Young diagram with  
the integer $\nu _i$ as the number of boxes in the $i$th column.
Thus, for example, $T(5,3,1)=\whT(3,2,2,1,1)$.

\item 
If $T=T(n_1, \cdots, n_r)=\whT(\nu_1, \cdots, \nu _p)$, then  
$n=n_1 + \cdots  +n _r$ and $n=\nu _1 \cdots +\nu _p$ are dual partitions to each other.

\item
When we say that $n=n_1 + \cdots + n_r$ is a partition of $n$, 
it is not assumed that the terms are arranged in  either decreasing or increasing 
order,  but it is assumed that each term is positive. 

\item 
When we know that the sequence in  a partition $n=n_1 + \cdots + n_r$  is put in  decreasing order,  
we use the term ``dual partition'' for $\whT(n_1, \cdots, n_r)$, identifying it with the Young diagram. 

\item 
When  $T=T(n_1, \cdots, n_r)$, sometimes $T$  is referred to as a sequence  
or a partition  in the obvious sense. 
\end{itemize}
\subsection{The linear hull  of a generic exponential matrix}
Until  Proposition~1 begins,  we assume that char $K$=0.  Let  $J \in \bM(n)$.  As is well known, 
the exponential of $J$ is the following:
$${\displaystyle \exp(J)=E+J+\frac{1}{2!}J^2+\frac{1}{3!}J^3+ \cdots    }$$
We are interested in the linear hull of the set $\{   \exp(xJ)  | x \in K \}$,   where  $J$  is  a nilpotent matrix.
Suppose for the moment that $J$ is a  single $n \times n$ Jordan 
block of a nilpotent matrix, so 
$$J=
\left( \begin{array}{rrrrr}
0          &    1     &  0    &   \cdots   & 0   \\
0          &    0     &  1    &   \cdots   & 0   \\
            &\ddots &        &                 &    \\
0          &    0     &  0    &   \cdots    & 1  \\
0          &    0     &  0    &   \cdots    &  0
\end{array}
\right). 
$$
Then 
$$\exp(xJ)=
\left( \begin{array}{lllll}
1    &    x     &  x^2 /2!  &   \cdots           & x^{n-1}/(n-1)!   \\
0    &    1     &  x         &   \ddots               & \vdots          \\
0    &    0     &  1         &   \ddots               & x^2/2!   \\
0    &    0     &  0         &   \ddots               & x         \\
0    &    0    &   0          &   0                         &  1
\end{array}
\right). 
$$
Thus the linear hull of $   \{ \exp(xJ) | x \in K       \}$
is the vector space spanned by $E, J, J^2, J^3, \cdots, J^{n-1}$.  
If $x$ is an indeterminate, we say 
$\exp(xJ)$ is a {\bf generic exponential } of $J$.
By an {\bf augmented exponential} we mean a matrix of type either   $(O|\exp (J))$ or 
its vertical version, 
$$\left( \begin{array}{c}
\exp (J)  \\ \hline
O
\end{array}
\right),
$$
where $O$ is a zero block of an arbitrary size as long as it fits $\exp(J)$.   
If $J$  is a single Jordan cell of a nilpotent matrix, 
the linear hull of an 
augmented generic exponential of $J$ is the vector space 
consisting of matrices  which are one of the following types.
\begin{equation}   \label{lower_cell}
\left( \begin{array}{lllllll}
0     &     0      &   0         & x_{0}    &    x_{1}      &  x_{2}         & x_{3}     \\
0     &     0      &   0         &   0         &  x_{0}        &  x_{1}         & x_{2}     \\
0     &     0      &  0          &    0         &      0           & x_{0}          & x_{1}    \\
0      &     0       &             0   &         0      &     0     &  0              & x_{0}      
\end{array}
\right) 
\end{equation}

\begin{equation}   \label{upper_cell}
\left( \begin{array}{llll}
x_{0}    &    x_{1}      &  x_{2}         & x_{3}     \\
0           &    x_{0}      &  x_{1}         & x_{2}     \\
0           &    0             &  x_{0}         & x_{1}     \\
0          &     0             &  0                & x_{0}     \\
0          &     0             &  0                & 0            \\
0          &     0             &  0                & 0      
\end{array}
\right) 
\end{equation}

\subsection{The commutator algebra of a nilpotent matrix as a set}    

For $J  \in \bM(n)$ we denote by $\C(J)$ the  commutator algebra of $J$, namely 
$$\C(J)=\{X \in \bM(n) | XJ=JX\}.$$
Note that $\C(J)$ is an associative algebra with identity.  

Let  $T=T(n_1, n_2, \cdots, n_r)$ be a Young diagram of size $n$ numbered horizontally  
and let $J \in \bM(n)$ be the Jordan canonical matrix of partition $(n_1, n_2, \cdots, n_r)$, so $J$ 
is  the  matrix defined by the equation~(\ref{jordan_matrix}).  Let $J_i$ be the $i$th diagonal block of $J$,  
namely,  $J_i$ is the  Jordan cell of size $n_i$. 
For the moment let us write  $\expo(J_i)$ for  $\exp(J_i)$ augmented by $O$. 
Recall that $\expo(J_i)$ is determined by its size.
(For definition see Section~2.2.)  
 Introduce letters $x_{ij}$ as many as the number of  pairs 
 $(i,j)$ for $1 \leq i,j \leq r$,
  and define, for each pair $(i,j)$,  the matrix    $X_{ij}$   of  
size $n_i \times n_j$ 
as follows: 
\begin{equation} \label{def_of_X_ij}
X_{ij}=\left\{
\begin{array}{l}
  \expo(x_{ij}J_j)  \mbox{ if $i \leq j$, } \\
  \expo(x_{ij}J_i)  \mbox{ if $i > j$.} 
\end{array}
\right.
\end{equation}
Put 
\begin{equation}\label{def_of_X} 
X=(X_{ij})
\end{equation} 
by which we mean the $n \times n$ block matrix with blocks $X_{ij}$ defined above.  
The following lemma was proved by   Turnbull and Aitken~\cite{hTaA32} and also by Gantmacher~\cite{fG59}.  
The proof here is due to Basili~\cite{rB03}.
\begin{lemma}[Gantmacher\cite{fG59}, Turnbull and Aitken\cite{hTaA32}] \label{description_of_C}  
The set $\C(J)$ coincides with the linear hull of the set   
$\{ X \in \bM(n) | x_{ij} \in K  \}$ as defined in $(\ref{def_of_X_ij})$  and  $(\ref{def_of_X})$.
\end{lemma}
{\em Proof.}
Suppose  $M$ is an $n \times n$ matrix.  
Let $M=(M_{pq})$ be the block decomposition  shown in the  picture below.
\begin{equation}  \label{general_block_matrix}
M=
\begin{array}{ccccc}   \\

                                                              \cline{1-4}
\svline{M_{11}}   &  \svline{M_{12}}   & \svline{\cdots}  & \svline{M_{1r}}  & \} n_1   \\   \cline{1-4}
\svline{M_{21}}   &  \svline{M_{22}}   & \svline{\cdots}  & \svline{M_{2r}}   & \} n_2   \\    \cline{1-4}
\svline{\vdots}   &  \svline{\vdots}   & \svline{\ddots}  & \svline{\vdots}   &  \vdots   \\    \cline{1-4}
\svline{M_{r1}}   & \svline{M_{r2}}      & \svline{\cdots}  & \svline{M_{rr}}   & \} n_r    \\    \cline{1-4}
\underbrace{}_{n_1}  &\underbrace{}_{n_2} & \ldots & \underbrace{}_{n_r}
\end{array}                    
\end{equation}
Then one sees that 
the condition $MJ=JM$ means that $M_{pq}J_q=J_pM_{pq}$ for all $1 \leq p,q  \leq r $. 
Thus the assertion follows from Lemma~\ref{commute} below.
\begin{lemma}   \label{commute}  
Let $Z=(z_{ij})$ be a $p \times q$ matrix.  
Let   $J_1$ and $J_2$   be the Jordan cells of sizes $p$ and  $q$ 
respectively.   
\begin{enumerate}
\item[$(i)$] If $p \geq q$, 
then $J_1Z=ZJ_2$ implies that 
$z_{21}=z_{31}=\cdots = z_{p1}=0$ and 
$z_{ij}=z_{(i+1)(j+1)}$ for all $i,j$ such that 
$1 \leq i \leq p-1$ and $1 \leq j \leq q-1$.  
\item[$(ii)$]  If $p < q$, 
then $J_1Z=ZJ_2$ implies that 
$z_{p1}=z_{p2}=\cdots = z_{p(q-1)}=0$ and 
$z_{ij}=z_{(i+1)(j+1)}$ for all $i,j$ such that 
$1 \leq i \leq p-1$ and $1 \leq j \leq q-1$.  
\end{enumerate}
\end{lemma}
{\em Proof } is straightforward.

\begin{remark}  \label{small_remark}   
{\rm  
\begin{enumerate}
\item[$(i)$]   
A generic matrix $X$ in  $\C(J)$ decomposes as $X=(X_{ij})$,  where  each block  $X_{ij}$   is of the type 
shown in  (\ref{lower_cell})  for $i \geq j$ and  in   (\ref{upper_cell})  for $i < j$.
Entries in the last columns in the blocks $X_{ij}$ for $i \geq j$ 
and those in the first rows in $X_{ij}$ for $i < j$ are algebraically independent. 
\item[$(ii)$] 
As before let $(n_1, \cdots, n_r)$ be the partition of the nilpotent matrix $J$ put in the Jordan canonical form.  
Put $A=K[z]/(z^{n_1})$ and $V =\bigoplus _{i=1} ^{r} K[z]/(z^{n_i})$.  Regard  $V$   as   an  $A$-module.  
Then we may use monomials as a basis of $V$ so that $J$ is the matrix for the multiplication map 
\[
\times  z : V \ra V. 
\] 

It is not difficult to see that $\C(J)$ coincides with $\End_A(V)$.  Details are left to the reader.  
In this paper this is used only to indicate another proof of Theorem~\ref{main_theorem}. 
(See the last paragraph preceding  Theorem~\ref{main_theorem}.)   
\end{enumerate}
}
\end{remark}

\begin{example}   \label{example_332}  
{\rm  
 Let $T=T(3,3,2)$.  Then a general element of $\C(J)$ is of the form 
\begin{equation}
\left( \begin{array}{ccc|ccc|cc}
a   &    a'     &  a''       & b  &  b'  &  b''  & c  & c'      \\
0   &    a      &  a'       & 0   &  b  &  b'   & 0  & c     \\
0   &    0      &  a        & 0   &  0  &  b   &  0 & 0     \\  \hline
d   &    d'     &  d''       & e  &  e'  &  e''  & f  & f'      \\
0   &    d      &  d'       & 0   &  e  &  e'   & 0 & f     \\
0   &    0      &  d        & 0   &  0  &  e   &  0 & 0     \\  \hline
0   &    g     &  g'       & 0  & h &h'  & i  & i'      \\
0   &    0     &  g        & 0  & 0 &h  & 0  & i     
\end{array}
\right) 
\end{equation}
}
\end{example}

\begin{example}  
{\rm 
 Let $T=T(4,2,1)$.  Then a general element of $\C(J)$ is of the form 

\begin{equation}
\left( \begin{array}{cccc|cc|c}
a   &    a'     &  a''     & a'''  &  b  &  b'    & c        \\
0   &    a      &  a'     & a''    &  0  &  b    & 0       \\
0   &    0      &  a      & a'     &  0  &  0   &  0      \\ 
0   &    0      &  0      & a     &  0   &  0   &  0      \\   \hline
0   &    0     &  d      & d'     &  e  &  e'   & f    \\   
0   &    0     &  0      & d      &  0  &  e   &  0    \\  \hline
0   &    0     &  0     & g      &  0  &   h  &  i      
\end{array}
\right) 
\end{equation}
}
\end{example}

\section{The structure of the commutator algebra of a nilpotent matrix}    
To say anything about the structure of the algebra $\C(J)$ we need to look at the 
Young diagram $T$ for $J$  more  in detail.  
Throughout this section we fix $T$ and $J$ as follows: 
\begin{enumerate}
\item
$T=T(n_1, \cdots, n_r)$ is a Young diagram of size $n$.
\item
$J$ is the Jordan canonical form of a nilpotent matrix of type $T$. 
\end{enumerate}
Let  $(f_1, f_2, \cdots, f_s)$  be the finest subsequence of $(n_1, \cdots, n_r)$ such that 
$f_1 > f_2 >  \cdots > f_s > 0$.  Then we can rewrite the same sequence 
$(n_1, \cdots, n_r)$ as

\begin{equation}   \label{grand_block_decomposition}   
(n_1, \cdots, n_r)=(\underbrace{f_1 , \cdots f_1}_{m_1}, \underbrace{f_2, \cdots, f_2}_{m_2}, \cdots ,
 \underbrace{f_s, \cdots, f_s}_{m_s}). 
\end{equation}

The integer $m_j$ is the multiplicity of $f_j$.  Let us call 
$m_1, \cdots, m_s$  the 
{\bf multiplicity sequence} of $T(n_1, \cdots, n_r)$.  Note that it  
gives us a partition of the number $r$ of rows of $T$,  namely,  $r=m_1+m_2+ \cdots + m_s$.

Recall that the Jacobson radical of a ring is defined to be the intersection of  
$${\rm ann}(M)$$  where ${\rm ann}(M)$ 
denotes the annihilator of the module $M$ and 
 $M$ runs  over all simple (right) modules.  
The Jacobson radical is a two sided  nilpotent ideal and 
if it is $0$ then the ring is said to be  semisimple.  (See e.g., early pages of Herstein~\cite{iH68}.)    
The following is a  known result. 
E.g., this can be implied by \cite[Theorem~3.5.2]{yDvK93} with the identification of $\C(J)$ with  
$\End_A(V)$, where $A=K[z]/(z^{n_1})$ as described in Remark~\ref{small_remark}~$(ii)$.       
We give a direct proof after Proposition~\ref{key_proposition} and Example~\ref{ex_322} 
where we show a number of properties of 
a generic matrix of $\C(J)$.  

\begin{theorem}   \label{main_theorem} 
Let $\C(J) \subset \bM(n)$ be the commutator algebra of $J$. 
Let $m_1, m_2, \cdots, m_s$ be the multiplicity sequence of $T$. 
Let $\rho$ be the Jacobson radical of $\C(J)$.    Then 
there is a surjective homomorphism 
$$\Phi : \C(J)  \ra  \bM(m_1) \times \bM(m_2) \times \cdots  \times \bM(m_s).$$
with $\ker \:  \Phi = \rho$. 
\end{theorem}

Proof is postponed to the end of Definition~\ref{def_diagonal_blocks}.  
The map $\Phi$ is defined in (\ref{big_phi}).  
     
Let $M \in  \bM(n)$.   Using the partition  $n=n_1+ \cdots +n_r$  
we decompose $M$  into blocks  as indicated in (\ref{general_block_matrix}). 
When $M \in \bM(n)$ is considered as a block matrix in this way for a Young diagram 
$T=T(n_1, \ldots, n_r)$, 
we write  $M=(x^{(kl)}_{ij})$, 
by which we mean that the element $x^{(kl)}_{ij}$ is the $(i,j)$-entry 
of the $(k,l)$-block of $M$.  
Note that $M$ has square diagonal blocks. 
For $M=(x^{(kl)}_{ij})$ we define the matrix   $\whM$  by 
$$\whM=(x^{(ij)}_{kl}).$$   The matrices 
$M$ and $\whM$ differ only by a certain permutation of rows and columns.  
To name the permutation explicitly, let $T_h$ and $T_v$ be the same Young diagram with 
the horizontal and vertical numberings respectively.
Let $\pi: T_h \ra T_v$ be the permutation of the integers $\{1,2 ,\cdots, n\}$ which,
 as in the picture below,  claims that 
if a box is numbered $i$ in $T_h$ then the same box is numbered $\pi(i)$ in $T_v$.  
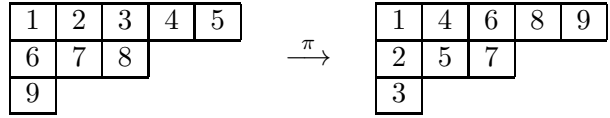
\begin{center}
\begin{tabular}{|ccccc}                                                                 \hline
1   &  \svline{2}   & \svline{3}  & \svline{4}  & \svline{5}     \\   \hline
6   &  \svline{7}   & \svline{8}  &                   &                       \\    \cline{1-3}
9   & \ssvline{}    &                   &                   &                        \\    \cline{1-1}
\end{tabular}                 \hspace{2ex}      $\stackrel{\pi}{\longrightarrow}$        \hspace{2ex}
\begin{tabular}{|ccccc}                                                                 \hline
1   &  \svline{4}   & \svline{6}  & \svline{8}  & \svline{9}     \\   \hline
2   &  \svline{5}   & \svline{7}  &                   &                       \\    \cline{1-3}
3   & \ssvline{}    &                   &                   &                        \\    \cline{1-1}
\end{tabular}                    
\end{center}
Let $P=(p_{ij})$ be the matrix defined by 
\begin{equation}   \label{permutation_matrix}
p_{ij}=\left\{\begin{array}{l}
1 \ \mbox{  if $j=\pi(i)$,}
\\
0 \ \mbox{  otherwise.}
\end{array}
\right.
\end{equation}
Then one sees easily that $\whM=P^{-1}MP$ for $M \in \C(J)$.  

Notice that the indices $i$ and $j$ of blocks  of $\whM=(a^{(ij)} _{kl})$ run over 1 through $n_1$, 
since $n_1$ is the size of the biggest block of $M=(a^{(kl)}_{ij})$.   
Let $n=\nu_1+ \cdots + \nu _p$ be the dual partition of 
$n=n_1+ \cdots + n_r$.  (So $\nu _i$ is the number of boxes of 
the $i$th  column of $T$.)
Then one sees with a little contemplation that 
the size of the $(i,j)$-block of $\whM$ is $\nu_i \times \nu_j$.  
We are interested in the diagonal blocks of $\whM$.  So we introduce a definition.

\begin{definition}  \label{en_sub_i}  
{\rm 
For a Young diagram $T=T(n_1, \cdots, n_r)$ and a block matrix $M=(a_{ij}^{(kl)})$ as in $($\ref{general_block_matrix}$)$ above, 
we define 
$$N_i$$
to be the $i$th diagonal block of $\whM$ 
for $ i=1,2,\cdots, p$.  Note that $N_i$ is a square matrix of size $\nu _i$, where $(\nu_1, \cdots, \nu_p)$ 
is the dual partition of $T=T(n_1, \cdots, n_r)$. 
}
\end{definition}

Recall that   $J$ is the  $n \times n$  Jordan matrix with partition $T=T(n_1, \cdots, n_r)$.  
\begin{proposition}  \label{key_proposition} 
Using the notation above  we have:
\begin{enumerate}
\item[$(i)$]
For $M \in \C(J)$, the matrix 
$\whM$ is  block upper triangular.  Namely, if $i > j$ the 
$(i,j)$-block of $\whM$ is O.

\item[$(ii)$]
$\whJ$ is the Jordan second canonical form of $T$.

\item[$(iii)$]
Let $\C(\whJ) \subset  \bM(n) $ be the commutator algebra of $\whJ$. 
Then the map $\C(J) \ra \C(\whJ)$ defined by    
\[
M \mapsto \whM
\] is 
a natural isomorphism of algebras. 

\item[$(iv)$]
Let $T'=T(n_1-1, n_2-1, \cdots, n_r-1)$, 
and let $J'$ be the Jordan canonical matrix with partition $T'$. 
Then 
$\whM$ decomposes as 
$$\whM=
\begin{array}{|c|c|}   \hline  
N_1   &  \ast   \\  \hline  
O     &  \whM'  \\  \hline 
\end{array}, 
$$  
where $M' \in \C(J')$. 

\item[$(v)$]
If $\nu_1=\nu_2$,  then $N_1=N_2$ for every $M \in \C(J)$.  
{\rm $($For the definition of $N_i$ see 
Definition~\ref{en_sub_i}.$)$}

\item[$(vi)$]
If  $\nu_1 > \nu_2$,  then $N_1$ decomposes as follows:

\begin{equation}   \label{decomposition_of_M_hat}
N_1=
\begin{array}{ccl}   \\

                                                              \cline{1-2}
\svline{N_{2}}  & \svline{G'}  & \} r-m_s       \\   \cline{1-2}
\svline{O}          &  \svline{G}   & \} m_s  \\   \cline{1-2}
\underbrace{}_{r-m_s}  & \underbrace{}_{m_s}
\end{array}                    
\end{equation}
Furthermore all entries of  $G$ are algebraically independent of 
any other entry of $M$ if $M$ is generic in $\C(J)$.  
\item[$(vii)$]
$N_1$ decomposes as:
\begin{equation}  \label{diagonal_N_1}
N_1=    
\begin{array}{cccccc}                                                                                                                                 \\   \cline{1-5}
\svline{G_{1}}           &  \svline{\ast}         & \svline{\ast}   &\svline{\ast}& \svline{\ast}      & \} m_1    \\   \cline{1-5}
\svline{O}                  &  \svline{G_{2}}       & \svline{\ast}   &\svline{\ast}& \svline{\ast}      & \} m_2    \\  \cline{1-5}
\svline{O}                  &  \svline{O}             & \svline{G_3}    &\svline{\ast}& \svline{\ast}     & \} m_3     \\  \cline{1-5}
\svline{\vdots}          &  \svline{\vdots}     & \svline{\vdots}  &  \svline{\ddots} & \svline{\vdots}     &  \vdots    \\   \cline{1-5}
\svline{O}                  & \svline{O}               & \svline{O}        &\svline{\cdots}& \svline{G_{s}}          & \} m_s \\  \cline{1-5}
\underbrace{}_{m_1} &\underbrace{}_{m_2} & \underbrace{}_{m_3} & \ldots  & \underbrace{}_{m_s}
\end{array}                    
\end{equation}
\end{enumerate}
\end{proposition}
{\em  Proof.}
$(i)$  Suppose $j>i$.  Let $x^{(ij)}_{kl}$  be the  $(k,l)$-entry 
of the $(i,j)$-block of $\whM$.  Originally it is the   
$(i,j)$-entry of $(k,l)$-block of $M$.
In Lemma~\ref{description_of_C} we showed that   
each block  is ``upper triangular.''   Hence  $x^{(ij)}_{kl}=0$ for any  $k,l$.  \newline
$(ii)$  Left to the reader.   \newline
$(iii)$  With $P$ as defined by       
(\ref{permutation_matrix}), we have $\whM=P^{-1}MP$. Hence the assertion follows.   \newline
$(iv)$  This is not difficult to see.  Details are left to the reader.  \newline  
$(v)$  That $\nu _1=\nu _2$ means that each block 
of $M$ has at least two rows and two columns.  Hence every block 
of $M$ has the $(2, 2)$  entry and moreover 
$x^{(kl)}_{11}=x^{(kl)}_{22}$ for all $(k,l)$.  This shows that 
$N_1=N_2$.   \newline
$(vi)$  That $\nu _1 > \nu _2$ means that $\nu _1 - \nu _2=m_s$ and that 
$f_s=1$.  So there are $m_s ^2$ blocks of size $1 \times 1$ in $M$. 
If we write $M$ with variable entries, then  these blocks  consist of algebraically independent elements by
Remark~\ref{small_remark}.  This  implies that the entries of $G_s$ consist of algebraically independent elements and  
are algebraically independent of any other entries of $M$.  
\newline
$(vii)$
In view of $(iv)$ the assertion follows inductively  from $(v)$ and $(vi)$.

\begin{example}   \label{ex_322} 
{\rm 
Let T=T(3,2,2).   Then a general element of $\C(J)$ is:

\begin{equation}
M=\left( \begin{array}{ccc|cc|cc}
a   &    a'     &  a''      &  b  &  b'    & c  &c'      \\
0   &    a      &  a'      &  0  &  b    & 0   &c    \\
0   &    0      &  a       &  0  &  0   &  0  &0    \\  \hline
0   &    d     &  d'       &  e  &  e'   & f    & f' \\   
0   &    0     &  d        &  0  &  e   &  0  & f  \\  \hline
0   &    g     &  g'       &  h  &   h'  &  i    & i' \\
0   &    0     & g        &  0  &   h   & 0    &  i  
\end{array}
\right) 
\end{equation}
By a  certain permutation it becomes $\whM$ as follows:
\begin{equation}
\whM=\left( \begin{array}{ccc|ccc|c}
a   &    b     &  c       & a'    &  b'  &  c'    & a''     \\
0   &    e     &  f       & d'    &  e'  &  f'    & d'      \\
0   &    h     &  i        & g     &  h'  &  i'   &  g'     \\   \hline
0   &    0     &  0       & a     &  b   &  c   &  d'      \\   
0   &    0     &  0      &  0     &  e  &    f   & d       \\   
0   &    0     &  0      &  0     &  h   &  i    &  g      \\  \hline
0   &    0     &  0      &  0     &  0  &   0  &   a      
\end{array}
\right) 
\end{equation}
Note that the block decomposition of $\whM$ gives us the dual partition $7=3+3+1$.
The first diagonal block of $\whM$ is the  matrix
$N_1=\left( \begin{array}{c|cc}
a   &    b     &  c       \\ \hline
0   &    e     &  f       \\ 
0   &    h     &  i       
\end{array}
\right)$.  The second $N_2$ is identical with $N_1$.  The third $N_3=(a)$. 
}
\end{example}
\begin{definition}   \label{def_diagonal_blocks}  
{\rm 
With the notation made in Definition~\ref{en_sub_i}, we call the sequence of matrices 
$$(N_1, N_2, \cdots, N_p)$$
the {\bf coarse diagonal blocks} of $\whM$.  
We call the sequence $(G_1, \cdots, G_s)$ the diagonal blocks of 
$N_1$.  We apply the  term analogously  to all $N_i$.  Hence the diagonal blocks of 
$N_i$ are $(G_1, G_2, \cdots, G_t)$ with a 
certain $t$ depending on $i$. 
By the {\bf  fine diagonal blocks} of $\whM$ we mean the totality of the diagonal blocks:
\begin{equation}  \label{coarse_diagonal_block}
(\mbox{\rm diag}(N_1), \mbox{\rm diag}(N_2), \cdots ,  {\rm diag}(N_p))  
\end{equation}
}
\end{definition}

{\em Proof}\ of Theorem~\ref{main_theorem}.
Let $M \in \C(J)$.  Let $N_1, N_2, \cdots, N_p$ be the coarse diagonal blocks of 
$\whM$, and $G_1,  \ldots, G_s$ be the diagonal blocks of $N_1$.
Define the map
\begin{equation}   \label{big_phi}
\Phi: \C(J)  \ra \bM(m_1) \times \cdots \times \bM(m_s) 
\end{equation}
by 
$\Phi(M)=(G_1, \cdots, G_s)$.   
By Proposition~\ref{key_proposition} $(vi)$ and $(vii)$  
it is easy to see that $\Phi$ is surjective and also that the kernel of $\Phi$ 
is nilpotent. 
Since $\bM(m_1) \times \cdots \times \bM(m_s)$ is semisimple, the proof  of Theorem~\ref{main_theorem} 
is complete.

\begin{remark}   
{\rm 
The matrices $G_{i}$ defined above are precisely the same as $\bar{A}_{\alpha \beta}$ with $\alpha=\beta=i$ 
in Basili~\cite[p.60]{rB03}.   
Basili shows that $M$ is nilpotent if and only if each $\bar{A}_{ii}$ is nilpotent.  
This is particularly obvious after the proof of Theorem~\ref{main_theorem}.  Thus we recover Basili's identification 
of the nilpotent commutator with the inverse image under $\Phi$ of the locus where the $\bar{A}_{ii}$ or $G_{i}$ are 
nilpotent.
} 

\end{remark}

In the next proposition we would like to redefine $\Phi$ 
in a coordinate free manner.  
Let $V$ be a finite dimensional  vector space over $K$.  
We use the same letter $J$ as before to denote a nilpotent element of $\End(V)$.  
The notation $\C(J)$ is used in the same meaning as for the matrix.
Namely $\C(J)= \{ M \in \End(V) | MJ=JM \}$.  
Note that the subspaces  $\ker J^i$ and  $\im  J^i$ of $V$ are  $\C(J)$-modules for every 
integer $i$, and  so are their sums and intersections.  
Let $p$ be the least integer such that $J^p  =0$.   
(To avoid the trivial case we assume $p>1$.)
We have a descending chain of subspaces:
\[
V=\ker J^{p} + \im J \supset \ker J^{p-1} + \im J \supset \ker J^{p-2}+ \im J \supset \cdots \supset \ker J^{0}+\im J=\im J
\]
From among the sequence of successive quotients 
$$(\ker J^{p-i} + \im J)/(\ker J^{p-i-1}+\im J)$$  for $i=0,1, \cdots, p-1$,  
pick the non-zero vector spaces  and rename them  
\begin{equation}  \label{definition_of_U}
U_1, U_2, \cdots, U_s
\end{equation}
Note that $U_1=V/(\ker J^{p-1}+ \im J)$.  In fact everything may be regarded as 
a module over $K[J]$, which is a commutative local ring.   
So we may use Nakayama's Lemma to see $V \not = (\ker J^{p-1}+ \im J)$.

Likewise consider the ascending chain of subspaces:
$$\label{definition_of_W}
0=\im J^{p} \cap \ker J \subset \im J^{p-1} \cap \ker J \subset \im J^{p-2} \cap \ker J \subset \cdots \subset \im J \cap \ker J \subset \im J^{0} \cap \ker J = \ker J
$$
Let 
\begin{equation}
W_1, W_2, \cdots, W_{s'}
\end{equation}
be the subsequence of non-zero terms of the successive quotients  
$$(\im J^{p-i-1} \cap \ker J)/(\im J^{p-i} \cap \ker J)   \mbox{ for } i=0,1, \cdots , p-1$$  
(We note  $W_1=\im J^{p-1} \cap \ker J$.)
The spaces  $U_i$ and $W_i$ are  $\C(J)$-modules.  Hence the  module structure induces  the 
algebra homomorphisms $\phi _i: \C(J) \ra \End(U_i)$ for $i=1,2, \cdots, s$.
Define the algebra homomorphism  
\begin{equation} \label{small_phi}
\phi : \C(J) \ra \End(U_1) \times \End(U_2) \times \cdots \times \End(U_s)
\end{equation}
by the concatenation $\phi= (\phi _1, \cdots, \phi _s)$ of all $\phi _i : \C(J) \ra \End(U_i)$.
Similarly define the algebra homomorphism  
\begin{equation} \label{small_phi_2} 
{\phi}' : \C(J) \ra \End(W_1) \times \End(W_2) \times \cdots \times \End(W_{s'})
\end{equation}
by the concatenation of ${\phi _i}' : \C(J) \ra \End(W_i)$.   

Now we have
\begin{proposition} \label{central_module} 
\begin{enumerate}
\item[$(i)$]
$s=s'$ 
\item[$(ii)$] 
$\dim U_i =\dim W_i$ for $i=1,2, \cdots, s$  
\item[$(iii)$]
If we identify  $\End(U_i)$ and $\End(W_i)$ with a full matrix ring using suitable bases, then the homomorphism  
$\Phi$ defined in $(\ref{big_phi})$  
coincides with $\phi$ in $(\ref{small_phi})$  and with  ${\phi'}$ in $(\ref{small_phi_2})$.   
\end{enumerate}
\end{proposition}

{\em Proof.}   Let $n=\dim V$.  
Let  $B$ be a basis of $V$ in which $J$ is put in the Jordan first canonical form.   Once and for all 
we fix such a basis and we identify  $\End(V)$ and $\bM(n)$.
Suppose that the matrix  for $J$ decomposes into  Jordan cells as denoted by the  
Young diagram  $T=T(n_1, n_2, \cdots, n_r)$.
We  index the boxes of $T$  by 
the basis elements of $V$ themselves  in such a way that it satisfies the following condition:
\[
e, e' \in B \mbox{ and }  e'=Je  \Leftrightarrow \mbox{ the box $e'$ is next to and on the right of  $e$.}  
\] 
(cf.  Equation~(\ref{jordan_matrix}) of Section~2.1.)
Here the ``box $e$'' means that the box is indexed by $e$.  Henceforth 
we use the words ``box in $T$'' and a ``basis element'' in $B$   interchangeably.  

Let $(m_1, m_2, \cdots, m_{s''})$ be the multiplicity sequence of $T$ so we have, by (\ref{grand_block_decomposition}) 
in Section~3, 
$$   
(n_1, \cdots, n_r)=(\underbrace{f_1 , \cdots f_1}_{m_1}, \underbrace{f_2, \cdots, f_2}_{m_2}, \cdots ,
 \underbrace{f_{s''}, \cdots, f_{s''}}_{m_{s''}}).
$$
Let us call the set of boxes in $T$ a ``rectangle''  if it consists of all rows of the same length $f_i$ 
for some $i$.  Let $T_i$ be the $i$th rectangle.  
(So $T_i$ consists of $m_i \times f_i$ boxes.)
We may write $T=T_1 \sqcup T_2 \sqcup \cdots \sqcup T_{s''}$ as a disjoint union of 
Young subdiagrams aligned left. 
Let $e$ be a basis element of $V$ in the first column of $T$ and let $e'$ be a basis element at the end of a row of $T$. 
We have that 
$$\mbox{the box $e$ belongs to }\ T_j  
\Leftrightarrow  \left\{  
\begin{array}{l}
J^{k}e=0  \mbox{ if $k=f_j$,} \\  
J^{k}e \not = 0 \mbox{ if $k < f_{j}$,}
\end{array}
\right.   
$$  
and 
$$\mbox{the box $e'$ belongs to }\ T_j  
\Leftrightarrow  \left\{  
\begin{array}{l}
e' \in \; \im J^k  \mbox{ if $k=f_j-1$,} \\  
e' \not \in \; \im j^k  \mbox{ if $k > f_{j}-1$.}
\end{array}
\right. 
$$  
Thus the basis elements in the first column of $T_i$ span the space $U_i$ 
  and those in the last column span $W_i$.  
This shows  $(i)$ and $(ii)$.  
To prove $(iii)$, we have to show that for $M \in \C(J)$, both 
$\phi_i(M)$ and ${\phi}'_i(M)$, for each $i$,  are represented by the matrix $G_i$ defined in 
Proposition~\ref{key_proposition} with a suitable basis for $V$.  
Number the boxes of $T$ vertically.  
Then one sees that       
 $\phi_i(M)$, for $M \in \C(J)$,  corresponds to the 
$i$th diagonal block  $G_i$ of $N_1$ as was shown 
in  Proposition~\ref{key_proposition} $(vii)$.  
Also It is not difficult to see that  ${\phi}'_i(M) \in \End(W_i)$ 
corresponds to the last fine 
diagonal block in  $N_{f_{s+1-i}}$, which is also $G_i$.   
This completes the  proof. 
\begin{remark}  \label{central_simple_module} 
{\rm 
Let $J \in \End(V)$ be nilpotent and let  
$U_1, \ldots, U_s$ and  $W_1, \ldots, W_s$ be the vector spaces defined above. 
Suppose that $T(n_1, \cdots, n_r)$ is  the Jordan decomposition 
for $J$.      Let $(f_1, f_2, \cdots, f_s)$ be the finest subsequence 
of $(n_1, \cdots, n_r)$ such that $f_1 > f_2 > \cdots > f_s >0$.  
Then we have:
\begin{enumerate}
\item[$(i)$]
$U_i=(\ker J^{f_i} + \im J)/(\ker J^{f_i -1} + \im J) \mbox{ for  $i=1,2,\cdots, s$. }$ 
\item[$(ii)$]
$W_i=(\ker J \cap \im J^{f_i-1})/(\ker J \cap \im J^{f_i}) \mbox{ for  $i=1,2,\cdots, s$. }$ 
\end{enumerate}
These could have been the definitions of the vector spaces $U_i$ and $W_i$.  
The statement of Proposition~\ref{central_module}~$(iii)$ means that $\phi _i$ and ${\phi}' _i$ are equivalent 
in the sense that there are bijective linear maps $\psi _i: U_i \ra W_i$,  which make the following  
 diagrams commutative 
\def\haphi{\stackrel{\phi _{i}(X)}{\longrightarrow}}
\def\haphiprime{\stackrel{\phi '_{i}(X)}{\longrightarrow}}
\newcommand{\mapdown}[1]{\downarrow\rlap{$\vcenter{\hbox{$\scriptstyle #1\, $}}$ }}
\def\dapsi{{\mapdown{\psi_i}}}
\noindent 
$${\displaystyle 
\begin{array}{ccccccccc}
    U_i        &     \haphi     &  U_i        \\[1ex]
    \dapsi     &                &  \dapsi     \\[1ex]
    W_i        & \haphiprime    &  W_i        \\[1ex]
\end{array}
}
$$
for every $X \in \C(J)$.  
This shows that $U_i$ and $W_i$ are isomorphic as modules over $\C(J)$.  
If we take grading into account, we have the isomorphism 
\[U_i[1-f_i]\cong W_i,\]
since the degrees of the boxes of the first column of $T_i$ are different 
from those in the last column of $T_i$ by $f_i-1$ row-wise. 
This isomorphism can be proved more directly as follows. 
\def\raa{\stackrel{J^{i}}{\ra}}
\def\rab{\stackrel{J^{i-1}}{\ra}}
\def\daa{{\mapdown{J}}}
\def\dab{{\mapdown{J}}}

\noindent 
Consider the diagram.  (We write $0:J$ for $\ker J$.)
$${\displaystyle 
\begin{array}{ccccccccccccccc}
  &      & 0                   &       &   0                  &          &                        &      &    \\[1ex]
  &      & \da                 &       & \da                  &          &                        &      &    \\[1ex]
0 & \ra  &  0:J                & \ra   & 0:J                  & \ra      &  0                     &      &    \\[1ex]
  &      & \da                 &       & \da                  &          &  \da                   &      &    \\[1ex]
0 & \ra  & 0:J^{i}             & \ra   & 0:J^{i+1}            & \raa     &  0:J \cap \im J^{i}    & \ra  & 0  \\[1ex]
  &      & \daa                &       & \dab                 &          &  \da                   &      &    \\[1ex]
0 & \ra  & 0:J^{i-1}           & \ra   & 0:J^{i}              & \rab     &  0:J \cap \im J^{i-1}  & \ra  & 0  \\[1ex]
  &      & \da                 &       & \da                  &          &  \da                   &           \\[1ex]
0 & \ra  & 0:J^{i-1}/J(0:J^i)  & \ra   & 0:J^{i}/J(0:J^{i+1}) & \ra      &  X                     & \ra  & 0  \\[1ex]
  &      & \da                 &       & \da                  &          &  \da                   &      &    \\[1ex]
  &      & 0                   &       &   0                  &          &  0                     &      &    
\end{array}
}
$$

\noindent
The definition of the maps should be self-explanatory. 
Also note the isomorphisms
\[
\frac{0:J^i}{J(0:J^{i+1})} \cong \frac{0:J^{i}+\im J}{\im J}
\]
and the same for $i-1$ instead of $i$. 
Then the last horizontal exact sequence shows that $X$ is isomorphic to $U_j$ if $i=f_j$.  
But the last vertical 
exact sequence shows that $X \cong W_j$.  Thus we have $U_j \cong W_j$ for all $j$.  

}
\end{remark}

\begin{remark} \label{new_remark}  
{\rm 
\begin{enumerate}
\item
$U_1, \cdots, U_s$ are simple $\C(J)$-modules of different isomorphism types, and these exhaust all 
isomorphism types of simple  $\C(J)$-modules. 

\item
When $J$ is given as the multiplication map $\times z \in \End(A)$, $A$ an Artinian $K$-algebra, 
we term the modules $U_1, \cdots, U_s$ the ``central simple modules'' of the pair $(A,z)$,  and 
study them further in \cite{tHjW06} and \cite{tHjW07}.  
We say  ``central simple'' in the sense it has no proper  submodules     
 over the centralizer of $J$.

\item  
Let $A$ be an Artinian Gorenstein  $K$-algebra, not necessarily graded,    
and let $z \in A$ be any non-unit element.  
Let 
$(f_1^{m_1}, \cdots, f_s^{m_s})$  be the partition for 
the Jordan decomposition of the nilpotent element $\times z \in \End(A)$.
Then it is possible to define $U_1, \cdots, U_s, W_1, \cdots, W_s$ as 
in Remark~\ref{central_simple_module}.  
Since  $A$ is Artinian  Gorenstein, ${\rm Hom}_A({}\_, A)$ is an 
exact functor.   Hence we have the isomorphism  
\[
{\rm Hom} _A (U_i, A) \cong W_i, 
\mbox { and }  
{\rm Hom} _A (W_i, A) \cong U_i.   
\] 
\noindent So 
\[
{\rm Hom} _A (U_i, A) \cong U_i,  \mbox{ and 
}  {\rm Hom} _A (W_i, A) \cong W_i.  
\]
This explains the symmetry of the Hilbert function of $U_i$ shown in \cite[Proposition~4.6]{tHjW06}  and 
\cite[Proposition~5.3]{tHjW07}.  
\end{enumerate}
}
\end{remark}   
\begin{proposition} \label{conjugation}  
Let $T$, $J$ and $\C(J)$ be the same as above.  Given $M \in \C(J)$, 
 it is possible to put the fine diagonal blocks of $\whM$ in the Jordan 
 canonical form by conjugation without affecting the shape of $\whJ$. 
 \end{proposition}
 
{\em Proof.}  We have to find an invertible matrix $H$ such that 
$H^{-1}\whJ H=\whJ$  and 
the fine diagonal blocks  of $(H^{-1}\whM H)$  are  Jordan first canonical forms.   Let 
$N_1, \cdots, N_p$ be the diagonal blocks as in the proof of Theorem~\ref{main_theorem}.  
Let $G_1, \cdots, G_s$ be the diagonal blocks of $N_1$.  
Let $F_i \in \bM(m_i)$ be an invertible matrix such that 
$F^{-1}_iG_iF_i$ is the Jordan canonical form of $G_i$.
Let $H_1\in \bM(r)$ be the matrix which has $F_i$ as the $i$th diagonal 
block and O off diagonal.  Then $H_1$ puts the diagonal blocks of $N_1$ into    
the Jordan canonical forms.  In the same way define $H_i$ for $N_i$ for 
each $1 \leq i \leq p$.  Finally define $H \in \bM(n)$ so that it has $H_i$ as the 
$i$th diagonal block and O off diagonal. 
One sees easily that $H$ does not change the shape of $\whJ$,  so it has the 
desired property. 
\begin{proposition} \label{rank_of_big_matrix} 
Let $T=T(n_1, n_2, \cdots, n_r)$, $J$ and $\C(J)$ be the same as before.  
Let $m_1, \cdots, m_s$ be the multiplicity sequence of $T$ so that we have 
$$ 
T(n_1, \cdots, n_r)=T(\underbrace{f_1 , \cdots f_1}_{m_1}, \underbrace{f_2, \cdots, f_2}_{m_2}, \cdots ,
 \underbrace{f_s, \cdots, f_s}_{m_s})
$$
where $f_1, \cdots f_s$ is the descending subsequence of $n_1, \cdots, n_r$.
Let $M \in \C(J)$ and let $N_1$ be the first coarse diagonal block of  $\whM$.
Suppose that the  diagonal blocks of $N_1$ are $(G_1, \cdots, G_s)$ and 
moreover that all entries of $\whM$ are 0 except in the  diagonals $N_1, \cdots, N_p$.
Then we have:  
$$\mbox{\rm rank}\ (\whM + \whJ)  \geq  
\mbox{\rm rank }G_1^{f_1} + \mbox{\rm rank }G_2^{f_2} + \cdots + \mbox{\rm rank }G_s^{f_s} + \mbox{\rm rank }J.$$
Assume furthermore that  all entries of  $N_1$ are zero outside of  the
diagonal blocks $G_1, \cdots, G_s$.  Then we have: 
$$\mbox{\rm rank}\ (\whM + \whJ) =  
\mbox{\rm rank }G_1^{f_1} + \mbox{\rm rank }G_2^{f_2} + \cdots + \mbox{\rm rank }G_s^{f_s} + \mbox{\rm rank }J.$$
\end{proposition}

{\em Proof.} 
First we prove the second assertion. 
Assume $s=1$.
It means that $T=T(\underbrace{n_1, \cdots, n_1}_{r})$.  In this case the matrix $\whM+\whJ$ has the form:

\begin{equation}   \label{matrix2}
\begin{array}{ccccccc}                                                                                                                                                              \\   \cline{1-6}
\svline{N_1}          &  \svline{E}       &  \svline{O}           & \svline{\cdots}  &\svline{O}  & \svline{O}                & \} m      \\   \cline{1-6}
\svline{O}              &  \svline{N_2}   &  \svline{E}            & \svline{\cdots}  &\svline{O}  & \svline{O}                & \} m    \\  \cline{1-6}
\svline{O}              &  \svline{O}        &  \svline{N_3}        & \svline{\cdots}  &\svline{O}& \svline{O}                   & \} m   \\  \cline{1-6}
\svline{\vdots}      &  \svline{\vdots}  &  \svline{\vdots}  & \svline{\ddots}  & \svline{\ddots} & \svline{\vdots}  &  \vdots     \\    \cline{1-6}
\svline{O}              & \svline{O}           & \svline{O}            & \svline{\cdots}  &\svline{N_{p-1}} & \svline{E}      & \} m   \\    \cline{1-6}
\svline{O}              & \svline{O}            & \svline{O}             & \svline{\cdots}  &\svline{O} & \svline{N_p}               & \} m   \\    \cline{1-6}
\underbrace{}_{m} &\underbrace{}_{m}  &\underbrace{}_{m} & \ldots                &                  & \underbrace{}_{m}
\end{array}                    
\end{equation}

Here $G_1=N_1=N_2= \cdots =N_p$, and $p=n_1=f_1$ and $m=r$.   
We are going to apply basic row and column operations to this 
matrix so that the rank can be computed.  
Use row operations by block so the matrix becomes:

\begin{equation}
\begin{array}{ccccccc}                                                                                                                                                              \\   \cline{1-6}
\svline{G}          &  \svline{E}       &  \svline{O}           & \svline{\cdots}  &\svline{O}  & \svline{O}                & \} m      \\   \cline{1-6}
\svline{-G^2}     &  \svline{O}        &  \svline{E}            & \svline{\cdots}  &\svline{O}  & \svline{O}                & \} m    \\  \cline{1-6}
\svline{G^3}       &  \svline{O}        &  \svline{O}        & \svline{\cdots}  &\svline{O}& \svline{O}                   & \} m   \\  \cline{1-6}
\svline{\vdots}      &  \svline{\vdots}  &  \svline{\vdots}  & \svline{\ddots}  & \svline{\ddots} & \svline{\vdots}  &  \vdots     \\    \cline{1-6}
\svline{\pm G^{p-1}}    & \svline{O}           & \svline{O}            & \svline{\cdots}  &\svline{O} & \svline{E}      & \} m   \\    \cline{1-6}
\svline{\mp G^p}          & \svline{O}            & \svline{O}             & \svline{\cdots}  &\svline{O} & \svline{O}               & \} m   \\    \cline{1-6}
\underbrace{}_{m} &\underbrace{}_{m}  &\underbrace{}_{m} & \ldots                &                  & \underbrace{}_{m}
\end{array}                    
\end{equation}
(We have put $G=G_1$.)  Now use column operations to kill all the matrices in the blocks of the first column,  except the 
block at the bottom. Thus the rank is equal to $\rank \, G^{p} + m(f_1-1)$ and the assertion 
of the proposition is proved for the case $s=1$. 
Now assume $s >1$.   Write $T$ as a disjoint union of rectangles, $T=T_1 \sqcup T_2 \sqcup \cdots  \sqcup T_s$.  
(It means that    $T_i=T(\underbrace{f_i, \cdots, f_i}_{m_i})$ for each $i$.)
We proceed to the general case by induction on $s$. 
As before let $N_1, N_2, \cdots, N_p$ be the coarse diagonal blocks of $\whM$ and let 
$I_1, I_2, \cdots, I_{p-1}$ be the above diagonal blocks for $\whJ$.  (These are 
described as matrices (3) and (4) in Section~2.2.)
Recall that $I_i$ is a matrix of size $\nu_i \times \nu_{i+1}$.  
If $\nu_i=\nu_{i+1}$, then $I_i$ is the identity matrix of that size 
and if $\nu_i > \nu_{i+1}$ then  it is the identity of size $\nu_{i+1}$ 
augmented by a zero block from below. 
Also recall that the diagonal blocks of $N_1$ are $G_1, \cdots, G_s$ as shown in the figure~(\ref{diagonal_N_1}). 
By Proposition~\ref{key_proposition},  $G_1$ is contained in $N_i$ for every $i$ as the first diagonal block.
Let us write $m:=m_1$ for the size of $G_1$.   Let $E^{(i)}$ be the submatrix of $I_i$ consisting 
of  the first $m$ rows and  $m$ columns.   ($E^{(i)}$ is nothing but the identity matrix of size $m$.)
Notice that the $G_1$ appears $p$ times as 
a fine diagonal block of $\whM$ (one time in every $N_i$), and that except the first 
there is  $E^{(i)}$  somewhere above $G_1$ in the same column block.  
Now making row operations using $E^{(i)}$ it is possible to kill all $G_1$ of the same column block.
As in the case $s=1$ we are left with consecutive powers of $G_1$ in the first column block. 
Now make column operations, using various $E^{(i)}$ to annihilate everything on the same row block. 
Then as a result all $G_1$ disappear except the single $G^p$ at the bottom of the 
first column block.  Moreover from every $N_i$ the first row block disappears.  
Now for the rows and columns that were not involved by the above procedure we may 
apply the induction hypothesis for  
$T(\underbrace{f_2, \cdots f_2}_{m_2}  \cdots  \underbrace{f_s, \cdots,f_s}_{m_s})$.
This completes the proof for the second statement.    The proof for the first inequality is proved similarly.

\section{Application to the theory of Artinian $K$-algebras}   

Before we state our main theorem we introduce some definitions. 
By dimension we will mean dimension as a $K$-vector space.  
\begin{definition} \label{def_cosperner} 

Let $A=\bigoplus _{i=0} ^c A_i$ be a graded Artinian $K$-algebra, where $A_0=K$ is a field and $A_c \not = 0$.
\begin{itemize}
\item
The algebra $A$ has the {\bf weak Lefschetz property} $(${\bf WLP}$)$  if there is a 
linear element
$y \in A_1$ such that the multiplication
$\times y : A_i \rightarrow A_{i+1}$ is either injective or surjective
for all $i=0,1, \cdots, c-1$.
A {\bf weak Lefschetz element} is a linear element $y$ with this property.
\item
The algebra $A$ has the {\bf strong Lefschetz property} $(${\bf SLP}$)$  if there is a 
 linear element
$y \in A_1$ such that the multiplication
$\times y^{c-2i} : A_i \rightarrow A_{c-i}$ is bijective 
for all $i=0,1, \cdots, [c/2]$.  
We call a linear element $y$ with this property a {\bf strong  Lefschetz element}. 

\item 
An element   $y \in A_1$  is a {\bf general linear form} if there is a non-empty Zariski open set $U \subset A_1$, 
such that $y$ has the same Jordan canonical form as every element $y' \in U$.  (A general linear form exists if 
$K$ is infinite.)  

\item
The {\bf Sperner number} of $A$ is
${\rm Max}\{ {\rm dim }A_i | i=0,1,2 \cdots , c \}$. 
\item
The {\bf CoSperner number} of $A$ is
${\displaystyle
\sum_{i=0}^{c-1} {\rm Min}\{
{\rm dim } A_i, {\rm dim } A_{i+1}
\}
}$.
\end{itemize}
\end{definition}

\begin{remark}  \label{remarks_on_SLP}  
{\rm
It is easy to see the following. 
\begin{enumerate}
\item[$(i)$]
For any $y \in A$, the rank of $\times y$ does not exceed the CoSperner number of $A$.
\item[$(ii)$]
For any $y \in A$, the dimension of $A/yA$ is no less than the Sperner number of $A$. 
\item[$(iii)$]
A linear form  $y \in A$ is a  weak  Lefschetz element for $A$  if and only if  
the rank of $\times y$  is equal to the CoSperner number of $A$.  
\item[$(iv)$]
If the Hilbert function of $A$ is unimodal, then we have 
$$\mbox{\rm Sperner} A+\mbox{\rm CoSperner} A= 
\mbox{\rm dim} A.$$
\item[$(v)$]
 If $A$ has the strong/weak Lefschetz property, then a general linear form is a strong/weak Lefschetz element. 
\item[$(vi)$]
 Suppose that  $y \in A$ is a linear form and $T=T(n_1, \cdots, n_r)$ is the  partition  for the 
nilpotent endomorphism $\times y \in \End(A)$.  Then $y$ is a weak Lefschetz element if and only if 
$r$, the number of the Jordan blocks of  $\times y$,  is equal to the Sperner number of $A$.  
Also $y$ is a strong Lefschetz element if and only if the dual partition $\whT(\nu _1, \cdots, \nu _p)$ of $T$  
is the one obtained from the unimodal Hilbert function of $A$.  (cf. \cite[Lemma~3.7]{tHjW07}.)    

\end{enumerate}
}
\end{remark}

We are going to apply Proposition~\ref{rank_of_big_matrix} to Artinian $K$-algebras to evaluate the rank of a 
general linear form.  
First we review some basic facts on the associated form ring of an Artinian algebra with respect to the 
principal ideal generated by a linear form.   This was motivated by the  necessity to prove 
Proposition~\ref{rank_of_deformation_2}.  

Let $A$ be a graded  Artinian $K$-algebra and let $z \in A_1$ be any linear 
form of $A$. 
Put 
$$\Gr_{(z)}(A)=A/(z)\oplus (z)/(z^2) \oplus (z^2)/(z^3) \oplus  \cdots \oplus  (z^{p-1})/(z^p).$$
Here $p$ is the least integer such that $z^p  =0$. 
As is well known $\Gr_{(z)}(A)$ is endowed with a commutative ring structure. 
The multiplication in $\Gr_{(z)}(A)$ is given by 
$$
(a+(z^{i+1}))(b+(z^{j+1}))=ab+(z^{i+j+1}), 
$$
where $a \in (z^i)$ and $b \in (z^j)$.
Note that $\Gr(A)$ inherits a grading from $A$  and 
in this sense $\Gr(A)$ and $A$ have the same Hilbert function. 
For  a non-zero element $a \in A$ there is 
$i$ such that $a \in  (z^i) \setminus (z^{i+1})$. 
In this case 
we write $a^{*} \in \Gr(A)$ for the natural image of 
$a$ in $(z^{i})/(z^{i+1})$. 

Let $\times : A \ra \End(A)$ be the regular representation of the algebra $A$.  
(So $\times a$  is the endomorphism  of $A$ defined by $\times a(b)=ab$ for $a,b \in A$.) 
Let $z$ be a linear form of $A$.  Since $A$ is Artinian,  $\times z$ is nilpotent. 
Let $T=T(n_1, \cdots, n_r)$ be the Young diagram for the Jordan 
canonical form of $\times z$.

One sees easily that the number $r$ of 
parts  of $T$  is equal to 
${\rm dim}\  {\rm ker}[\times z: A \ra A] ={\rm dim}A/(z)$, since each Jordan cell  of $\times z$ 
contributes 1 to the dimension of the kernel.   
The Young diagram  $T(n_1-1, n_2-1 , \cdots, n_r-1)$, with zero's deleted,  corresponds to 
the Jordan canonical form of the induced map  $\times {\overline z} \in \End(A/(0:z))$.  
Thus, inductively,   it follows that 
the dual of the partition $T(n_1, \cdots, n_r)$  is
$T(\nu_1, \nu_2, \cdots, \nu_p)$, with the integers
$\nu _i ={\rm dim}(z^{i-1})/(z^{i})$.

Let $\bB \subset A$ be a $K$-basis of $A$  in which 
$\times z$ is written as a Jordan canonical form.  
We identify the boxes of $T=T(n_1, \cdots, n_r)$  and the elements of $\bB$.  
With this identification a row of $T$  is a basis for a Jordan cell of $\times z$. 

Let $\bB_i=\bB \cap ( (z^{i-1})  \setminus (z^{i})  )$. It is easy to see that 
$\bB_i \sqcup \bB_{i+1} \sqcup \cdots \sqcup \bB_{p}$ 
is a $K$-basis for the ideal $(z^{i-1})$.
With the identification of $\bB$ and the boxes of $T$, the set $\bB_i$ corresponds to 
the boxes in the $i$th column of $T$.    

Now let $\bB^{\ast}$ be the natural image of $\bB$ in $\Gr(A)$, i.e, 
$\bB^{\ast} = \{b^{\ast}  \in \Gr(A) | b \in \bB\}$.   
Similarly let $\bB_i^{\ast}=\{b^{\ast} | b \in \bB_i\}$.  
One sees immediately that $\bB^{\ast}$ is a basis of $\Gr(A)$ in which 
the map $\times z^{\ast}$ is represented by a Jordan canonical form. 
It is also immediate to see that $\bB^{\ast}_i$ is a $K$-basis for $(z^{i-1})/(z^{i})$.   

Now we would like to prove  
\begin{proposition}  \label{application}   
We use the notation  above. 
\begin{enumerate} 
\item[$(i)$]
The linear maps 
$\times z \in \End(A)$ and $\times z^{\ast}\in \End(\Gr_{(z)}(A))$ have the 
same Jordan canonical form. 
\item[$(ii)$]
Let $y \in A$ be any element.   Then $\times y \in \C(\times z)$ and 
$\times y^{\ast} \in \C(\times z^{\ast})$. 

\item[$(iii)$]
Let $y \in A$ be any linear form. 
Let $P$ be the matrix for $\times y$ with the basis $\bB$ and similarly 
$Q$ the matrix for $\times y^{\ast}$ with the basis  $\bB^{\ast}$.
Then the coarse diagonal blocks of $\widehat{P}$ and those of $\widehat{Q}$ coincide.

\item[$(iv)$] The kernel of the multiplication map $z^{\ast}: \Gr_{(z)}(A) \ra \Gr_{(z)}(A)$  is 
given by 
$$\bigoplus ^{p} _{\alpha =1}\left((z^{i-1}) \cap (0:z)+ (z^{i})\right)/(z^{i})$$
\end{enumerate}
\end{proposition}

{\em Proof.} $(i)$ Consider the ideal of $\Gr(A)$ generated by a power of $z^{\ast}$.  
First note that $(z^{\ast})^{\alpha}= (z^{\alpha})^{\ast}$.   Now it is easy to see that
$(z^{\ast})^{\alpha}\Gr_{(z)}(A) \cong  (z^{\alpha})/(z^{\alpha +1})\oplus \cdots $, which 
implies that ${\rm rank} (\times z)^{\alpha}= {\rm rank} (\times z^{\ast})^{\alpha}$ for all $\alpha = 1,2, 3, \dots$
This shows that they have the same Jordan canonical form. 

$(ii)$ Trivial.

$(iii)$ 
A coarse diagonal block of $\widehat{P}$ is the matrix for the induced map 
$\times y :A \in \End{((z^{\alpha}})/(z^{\alpha +1}))$ with the basis 
$\bB_{\alpha}^{\ast}$.  Thus the assertion follows immediately. 

$(iv)$  Left to the reader.

Let $R=K[x_1, x_2, \cdots, x_d]$ be the polynomial ring and let $I \subset R$ be 
a homogeneous ideal such that $A=R/I$ is an Artinian $K$-algebra.  
Put $Z=x_d$.  For any homogeneous element $f \in R$ 
it is possible to write uniquely 
$$f= f_0 + f_1Z + f_2 Z^2 + \cdots + f_k Z^k$$
where  $f_i$  is a homogeneous polynomial in $K[x_1, \cdots, x_{n-1}]$.
Let $i$ be the least integer such that $f_i \not  = 0$.
In this case we will write ${\rm In}'(f)= f_iZ^{i}$. 
Furthermore we define 
${\rm In}'(I)$ to be the ideal of $R$ generated by the set 
$\{ {\rm In}' (f) \}$,  
where $f$ runs over homogeneous elements of $I$. 
It is well known that $R/{\rm In}'(I) \cong {\rm Gr}_{(z)}(A)$. 
(Here we have set $z= Z\  {\rm mod }\  I$)  Suppose that $d=2$.  Then one notices that 
${\rm In}'(I)$ coincides with the ideal generated by the initial terms of $I$ with respect to the 
reverse lexicographic order with $x_1 > x _2$.  The same notation, ${\rm In}'(I), {\rm  \ and  \  }\Gr_{(z)}(I)$  will be applied for a 
graded submodule $I$ of a finite colength in a free $R$-module.  

\begin{proposition} \label{rank_of_deformation_2}  
Let $K$ be an infinite  field, and  let 
 $V$ be a finite vector space over $K$.  Let $J \in \End(V)$  be nilpotent.  
 Choose a  basis of $\, V$ so  that we may identify $\End(V)=\bM(n)$,  where $n= \dim V$ and 
  $J$ is put in the Jordan first canonical form.  
Let  $\C(J) \subset \End(V)$ be the  commutator algebra of $J$. 
 Let $M \in  \End(V)=\bM(n)$ be nilpotent such that $ M \in \C(J)$.   Let $N_1, \ldots, N_p$ be 
the coarse diagonal blocks of $\whM$.  Let $M^{\dagger}$ be the matrix such that 
$\whM ^{\dagger}={\rm diag}(N_1, \cdots, N_p)$.
Then $M^{\dagger} \in \C(J)$.  Moreover we have 
$$
\mbox{\rm rank} (M^{\dagger}+ \lambda J) \leq  \mbox{\rm rank} (M + \lambda J) 
$$
for most  $\lambda \in K$.  
\end{proposition}

{\em Proof.}
It is easy to see that  $M^{\dagger} \in \C(J)$ so we omit the proof.  To prove the second assertion 
let $R=K[y,z]$ be the polynomial ring in two variables.  Define an algebra homomorphism $R \ra \End(V)$  
by $y \mapsto M$ and $z \mapsto J$.  Then we may regard $V$ as an  $R$ module with support in the maximal ideal $(y,z)$. 
(Note  that $V$ is not necessarily graded.) 
 Now $M$ is the matrix for the multiplication map $\times y : V \ra V$ and $J$ for $\times z: V \ra V$.  
Let $$\Gr_{(z)}(V)= V/zV\oplus zV/ z^2V \oplus \cdots \oplus z^{p-1}V/z^pV$$
The module $\Gr_{(z)}(V)$ has naturally  the structure of  $R$-module. 
One notices that $\times z  \in \End(\Gr_{(z)}(V))$
has the same Jordan canonical form as $J$.  Moreover the matrix for $\times y \in  \End((\Gr_{(z)}(V))$ 
is $M^{\dagger}$.  
Let $g \in R$ be a general linear form.   
Now  by \cite[Proposition~3.3]{tHjW06}  we have 

$$\dim V/g V \leq \dim \Gr_{(z)}(V)/g \Gr_{(z)}(V) $$

This proves the assertion as we may assume that $g=y + \lambda z$
for a sufficiently general $\lambda \in K$ .  

\begin{theorem}   \label{rank_of_general_element}   
Let $K$ be an infinite field and 
let  $A=\bigoplus A_i$ be a graded Artinian $K$-algebra and let $z \in A$ 
be  any  linear form.   
Suppose that the Jordan decomposition of the nilpotent element 
$$\times z \in \End(A)$$
is given by 
$$ 
T=T(n_1, \cdots, n_r)=T(\underbrace{f_1 , \cdots f_1}_{m_1}, \underbrace{f_2, \cdots, f_2}_{m_2}, \cdots ,
 \underbrace{f_s, \cdots, f_s}_{m_s}).  
$$
Let $y \in A$ be a linear form linearly independent of $z$.  
Let $J, M  \in \End(A)$ be the matrices  for  
$z, y$ with a basis of $A$ so that 
$J$ is in the Jordan second canonical form.   
Let $N_1, \cdots, N_r$ be the coarse diagonal blocks 
of $\whM \in \C(\whJ)$ $($as defined in 
{\rm Definition~\ref{def_diagonal_blocks}}$)$ and  
let 
 $G_1, \cdots, G_s$ be the diagonal blocks of $N_1$.  
Then we have: 
\begin{enumerate}
\item[$(i)$]
$\mbox{\rm rank }G_1^{f_1} + \mbox{\rm rank }G_2^{f_2} + \cdots + \mbox{\rm rank }G_s^{f_s} + \mbox{\rm rank }J 
\leq  \rank(M+ \lambda J) $ for most  $\lambda \in K$. 
\item[$(ii)$]
The equality  
\[
\mbox{\rm rank }G_1^{f_1} + \mbox{\rm rank }G_2^{f_2} + \cdots + \mbox{\rm rank }G_s^{f_s} + \mbox{\rm rank }J 
={\rm CoSperner}(A), 
\]
implies that  $y+ \lambda z$ is a weak Lefschetz element of $A$ for most $\lambda \in K$.  
\end{enumerate}
\end{theorem}
\def\B{{\bf B}}
{\em Proof.}
$(i)$
Put $G=\Gr_{(z)}(A)$  and let $z^{\ast}, y^{\ast} \in G$ be the initial forms of $z, y$ respectively.  
Recall that $\times z \in \End(A)$ and $\times z ^{\ast} \in \End(G)$ have the same Jordan canonical form. 
We may choose a basis $\B \subset A$ such that $J:=\times z$ is in Jordan canonical form  
as well as  $\times z^{\ast}$  with $\B ^{\ast}$. 
Let $M^{\dagger}$ be the matrix for $\times y^{\ast}$ with $B^{\ast}$.  
Since $y \not \in (z)$, 
the matrix  $\whM^{\dagger} \in \C(\whJ)$ consists of only diagonal blocks by 
the way the multiplication is defined in $G$.   
Moreover they are the same as those of $\whM$.  
Thus we have $$ \rank(M^{\dagger} + \lambda J)  \leq \rank(M + \lambda J) $$
for most of $\lambda \in K$ by  Proposition~\ref{rank_of_deformation_2}.    
Now the first inequality immediately follows from Proposition~\ref{rank_of_big_matrix}. 
$(ii)$ This follows form Proposition~\ref{remarks_on_SLP} $(iii)$. 

The following Theorem was proved in  \cite[Theorem~1.2]{tHjW06}.  The proof is essentially the same as 
that of Theorem~\ref{rank_of_general_element}~$(ii)$ above.  

\begin{theorem} \label{thm_from_jpaa} 
Let    $A$   be an Artinian Gorenstein $K$-algebra and $z \in A$ a linear  form. 
Let  $$U_1, \cdots, U_s$$
be the central simple modules  defined in $(\ref{definition_of_U})$  
for the nilpotent endomorphism $\times z \in \End(A)$.  Suppose that  
all $U_i$ have the strong Lefschetz property as $A$-modules.  
Then $A$ has the strong Lefschetz property. 
\end{theorem}

\begin{remark} \label{rem_on_csm} 
{\rm    
Let    $A$   be an Artinian  $K$-algebra and $z \in A$ a linear  form. 
The central simple modules of $(A,z)$ are defined to be the non-zero modules  
of the form 
$(0:z^f + (z))/ (0:z^{f-1}+(z))$.  
They are modules over the algebra $A/(z)$ and are determined by the 
Jordan canonical form of $\times z$.    
Let $G=\Gr_{(z)}(A)$  be the associated form ring. Then the endomorphims 
$\times z^{\ast} \in \End(G)$  
and  $\times z \in \End(A)$  
have the same Jordan canonical form.  
Thus the central simple modules of $(G, z^{\ast})$ can be  regarded as the same 
 modules over of $A/(z)$ with the identification $G/(z^{\ast})=A/(z)$.    
Suppose that $A$ is Gorenstein. Then, even though $G$ may not be Gorenstein, the strong 
Lefschetz property of the central simple modules of $(A,z)$ implies that 
$G$ has the strong Lefschetz property. 
(See \cite[Theorem~5.2]{tHjW07}.)   
}
\end{remark}  

\section{Examples}

In the following examples we show how Theorem~\ref{rank_of_general_element} can be used to 
compute the rank of a general linear form for $A$ and to prove the weak Lefschetz property  of  $A$.   
We proved in \cite{tHjMuNjW01}
that every Artinian complete intersection in codimension three over a field of  characteristic zero has the WLP.  
The method to prove it in Examples~\ref{ex_b3} and \ref{tama_sampia_2003} are different from the one used in 
\cite{tHjMuNjW01}. 
\begin{example}  \label{ex_b3}   
Assume that  $K$ is  a field of {\rm char}$\, K \neq 2$.  
Let   $R=K[x,y,z], I=(x^2+y^2+z^2, x^4+y^4+z^4, xyz)$, and $A=R/I$.  
$($We use the same letter $z$ for the image of $z$ in $A$.$)$  
Then $\times z \in \End(A)$ is 
represented by the partition 
\[
24=\underbrace{5+5+5+5}_{4}+\underbrace{1+1+1+1}_4
\]
The dual partition is 
\[
24=8+4+4+4+4
\]
The Young diagram is as follows:

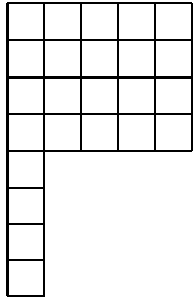
\begin{center}
\begin{equation}
\begin{array}{|ccccc}                                                                 \hline
{\ }    &  \svline{\ }   & \svline{\ }    &   \svline{\ }  &  \svline{\ }    \\   \hline
{}   &  \svline{}   & \svline{}     &   \svline{}  &  \svline{}   \\   \hline
{}  &  \svline{}   & \svline{}     &   \svline{}  &  \svline{}  \\   \hline
{}  &  \svline{}   & \svline{}      &   \svline{}  &  \svline{}   \\   \hline
{}   & \ssvline{}    &      & &                        \\    \cline{1-1}
{}   & \ssvline{}    &      & &                        \\    \cline{1-1}
{}   & \ssvline{}    &       & &                       \\    \cline{1-1}
{}   & \ssvline{}    &      & &                       \\    \cline{1-1}
\end{array}   
\hspace{8ex}
\end{equation}
\end{center}

The rank of a general linear form  of $A$ is $18$.   
$A$ has the strong Lefschetz property, but 
$z$ itself is not a strong Lefschetz element, since the rank of $\times z$ is 14.  

\end{example}

{\em Proof.}  
Since~char$\, K \neq 2$, we have that $A$ is Artinian. 
Consider the exact sequence
\[
0 \ra A/0:z \ra A \ra A/(z) \ra 0
\]
The first column of the Young diagram corresponds to $A/(z)$,  with the 
Hilbert function $1+2t+2t^2+2t^3+t^4$.  So the dimension is 8.
Put $B=K[x,y,z]/(x^2+y^2+z^2, xy, z^4)$.  
Then it is easy to see that 
$zA \cong A/(0:z) \cong K[x,y,z]/(x^2+y^2+z^2, xy, x^4+y^4+z^4) \cong B$.  
The ideal $zA$ corresponds to the diagram with the first column deleted.  
Now it is easy to compute $\dim B/(0:z^i)=4(4-i), \mbox{ for }  i=0,1,2,3$. 
Thus we have verified the partition for $\times z$ is 
$\whT(8,4,4,4,4)=T(5,5,5,5,1,1,1,1)$. 
Put $J= \times z$.  Then a general member of $\C(J)$ has coarse diagonal blocks 
$N_1, \cdots, N_5$  whose sizes are $(8,4,4,4,4)$ respectively  and $N_1$ has two diagonal 
blocks $G_1$ and $G_2$ of size $4$.  
Let $U_1$ and $U_2$ be the central simple  modules as defined in Remark~\ref{central_simple_module} and 
let $g$ be a general linear form of $A$.  
  Let $G_1$ be a  matrix for the induced map 
$\times g \in \End(U_1)$ and $G_2$ for $\times g \in \End(U_2)$.   
Notice that we have the exact sequence 
$$0 \ra U_2 \ra A/(z) \ra U_1 \ra 0,$$ where the first map sends 1 to $xy$.
Thus we have 
\[
\left\{\begin{array}{l} 
U_1 \cong K[x,y]/(x^2+y^2, xy),  \\ 
U_2 \cong K[x,y]/(x^2+y^2, xy)[-2]. 
\end{array}
\right. 
\] 
In the notation of Proposition~\ref{rank_of_big_matrix}, $s=2, f_1=5, f_2=1$ and 
$\rank G_1 ^{f_1} + \rank G_2 ^{f_2} + \rank J= 0+2+16=18$. 
Since the CoSperner number of $A$ is 18, this shows that  $A$ has  the weak Lefschetz property 
by Theorem~\ref{rank_of_general_element}.   By direct computation or using \cite{aI94} Theorem~2.9 or 
\cite{tHjMuNjW01}  Proposition~4.4, it follows that $U_1$ and $U_2$ have the SLP.  
Hence by Theorem~\ref{thm_from_jpaa}, $A$ has the SLP. 

\begin{example}  \label{tama_sampia_2003} 
{\rm 
Assume that $K$ is an infinite field of  $\mbox{char}\ K \neq 2$.  Let $R=K[x,y,z]$ and 
let $A=R/(x^4+y^4+z^4, xy^3+x^2z^2, y^3z)$.  Then, as is easily calculated,  $\times z \in \End(A)$ is represented by 
$$T=T(7,7,7,7,7,7,3,3,3,3,3,3,1,1,1,1)=\whT(16,12,12,6,6,6,6).$$
Thus using the notation of Proposition~\ref{rank_of_big_matrix}, $f_1=7, f_2=3, f_3=1$. 
Put $J=\times z$.  Let $M$ be a general member of $\C(J)$.      
The first coarse diagonal block $N_1$ of $M$ is of size $16$ and 
it consists of fine diagonal blocks $G_1, G_2, G_3$ of sizes 
$6,6,4$ respectively.  Let $U_1, U_2, U_3$ be the central simple modules  defined in Remark~\ref{central_simple_module}.  
Then, as with the previous example,  it is not difficult to see that $A$ is Artinian and that 
\[
\left\{
\begin{array}{l}
U_1 \cong K[x,y]/(x^2, y^3),       \\
U_2 \cong K[x,y]/(x^2, y^3)[-2],   \\
U_3 \cong K[y]/(y^4)[-3].          
\end{array}
\right.   
\]
Let $g \in A$  be a general linear form and let $G_i$ be the matrix for the 
induced map $\overline{g} \in  \End(U_i)$.  
Thus one sees that $G_1$ and $G_2$ are the nilpotent matrix  with Jordan decomposition $T(4,2)$ and 
$G_3$ with  $T(4)$.  
Thus $$\rank \: G_1 ^{f_1} + \rank \: G_2 ^{f_2} + \rank \: G_3 ^{f_3} + \rank \: J=0+ 1+ 3+ 48=52.$$  
This is equal to the CoSperner number of $A$. 
By Theorem~\ref{rank_of_general_element},  this shows that $A$ has the weak Lefschetz property.  
As in the previous Example, $z$ is not a weak Lefschetz element, but by  Theorem~\ref{thm_from_jpaa},  
$A$ has the strong Lefschetz property.
}
\end{example}

\begin{example}  \label{four_variable_yokohama} 
Let $R=K[w, x, y, z]$ be the polynomial ring, and   
put  $$A=R/(w^2, wx, x^3, xy, y^3, yz, z^3).$$
\end{example}

(We use the same letters $w,x, \cdots$ to denote their images in $A$.)
The Jordan decomposition of $J:=\times z \in \End(A)$ is represented by 
the partition $T=T(3,3,3,3,1,1,1,1)=\whT(8, 4, 4)$. 
For a general linear form $g \in A$, the matrix  $\times g \in \C(\whJ)$   has three coarse diagonal blocks of 
sizes $8,4,4$.  The first block $N_1$ has two diagonal blocks of 
size four each. 
One sees  that 
\[
\left\{ 
\begin{array}{l}
U_1 =A/(0:z^2 + (z)) \cong K[w,x,y,z]/(w^2, wx, x^3,  y, z),  \\
U_2= ((0:z) + (z))/(z) \cong K[w,x,y]/(w^2, x, y^2)[-1]. 
\end{array}
\right.
\]
(To see this,  notice that $(0:z) + (z))/(z)$ is a principal ideal of $A/(z)$ generated by $y$.)
Both $U_1$ and $U_2[1]$ have the Hilbert function $(1,2,1)$. 
Let $g \in A$ be a general linear form let  $G_i$  be the matrix for the induced  maps 
$\times \overline{g} \in \End(U_i)$.
Then $\rank (G_1) ^3 + \rank G_2 + \rank J= 0 + 2 + 8 = 10$, which is equal to the CoSperner number of $A$. 
Hence $A$ has the WLP.  In fact $A$ has the SLP, but Theorem~\ref{thm_from_jpaa} does not apply since $A$ is not 
Gorenstein. However, \cite[Theorem~5.2]{tHjW07} does apply.    

Alternatively we may let $w$ do the role of $z$.  
The Jordan decomposition for  $\times w$  is given by 
$T=T(\underbrace{2,2,2,2,2}_5,\underbrace{1,1,1,1,1,1}_{6})=\whT(11,5)$. 
We have   
\[
\left\{ 
\begin{array}{l}
U_1 =A/((0:w) + (w)) \cong K[w,x,y,z]/(w, x, y^3,  yz, z^3), \\
U_2= ((0:w) + (w))/(w) \cong K[x,y,z] / (x^2, y, z^3)[-1]. 
\end{array}
\right.
\]
$U_1$ has the Hilbert function $(1,2,2)$, and $U_2$ $(0, 1,2,2,1)$. 
(Notice that  $(1,2,2,0,0)+(0,1,2,2,1)=(1,3,4,2,1)$  is the 
Hilbert function of $A/(z)$.) 
Let $g \in A$ be a general linear form and let $G_1$ and $G_2$ be matrices for $\times \overline{g} \in \End(U_1)$   
and  $\times \overline{g} \in \End(U_2)$ 
respectively. Then 
$\rank (G_1) ^2 + \rank G_2 + \rank (J)= 1+ 4+ 5=10$.  This also shows that $A$ has the WLP\@.  
Note that we can use neither Theorem~\ref{thm_from_jpaa} nor \cite[Theorem~5.2]{tHjW07} to prove 
$A$ has the SLP, since $A$ is not Gorenstein and since $U_1$ does not have a symmetric Hilbert function. 

\begin{example}  \label{yokohama} 
Let $K$ be a field of characteristic 0. Let $R=K[x_1, x_2, \cdots, x_n]$ be the polynomial ring with $n \geq 2$.    
Let $$I=(x_1^2, x_1x_2, x_2^3, x_2x_3, x_3^3, \cdots, x_{n-1}x_{n}, x_n^3)$$
Then $R/I$ has the SLP.  
\end{example}

To prove this we  first consider  a similar but simpler example as follows.
\begin{example}   \label{simplified_yokohama}  
Let $K$ be a field of characteristic 0.  Let $n > 1$  be an integer. 
\begin{enumerate}   
\item[$(i)$]
$K[x_1, x_2, \cdots, x_n]/(x_1^2, x_2^2,  \cdots, x_{n-1}^2, x_{n-1}x_{n}, x_n^3)$ has the SLP.  
\item[$(ii)$]
$A=K[x_0, x_1, \cdots, x_n]/(x_0 ^{\alpha}, x_1^2, x_2 ^2, \cdots, x_{n-1}^2, x_{n-1}x_{n}, x_n^3 )$ has the WLP for any positive
 \linebreak[3] integer $\alpha$. 
\end{enumerate}
\end{example}

By \cite[Proposition~18]{tHjW03},  $(i)$ follows from $(ii)$.  
To prove $(ii)$ we would like to use  Theorem~\ref{rank_of_general_element} so that the same proof works for 
Example~\ref{yokohama} also. 
Put $z=x_n$. Note that $\dim A = 2^n\alpha$.   
Furthermore note that $\dim A/(z)= 2^{n-1}\alpha$ and $\dim (z)/(z^2) = \dim (z^2)/(z^3) = 2^{n-2}\alpha$.  
This shows  that the Jordan canonical form for $\times z \in \End(A)$ is given by the partition:
\[
T=T(\underbrace{3,3,\cdots ,3}_{2^{n-2}\alpha}, 
\underbrace{1,1,\cdots , 1}_{2^{n-2}\alpha})=\whT(\underbrace{2^{n-1}\alpha}_{1},\underbrace{2^{n-2}\alpha ,2^{n-2}\alpha }_{2})
\]

\noindent The Young diagram is shown in the picture below: 
\begin{equation}   \label{new}
\begin{array}{cccl}   \\
                                                              \cline{1-3}
\svline{\vdots}    & \svline{\ \ \vdots \ \ }    &    \svline{\ \  \vdots \ \ }     & \} {2^{n-2}\alpha}   \\   \cline{1-3}
\svline{\vdots}          &                            &                        & \} {2^{n-2}\alpha}        \\   \cline{1-1}
\underbrace{}_{1}  &    \multicolumn{2}{c}{  \underbrace{\ \ \ \ \ \ \ \ }_{2}   }
\end{array}                    
\end{equation}

Now we see that the first coarse diagonal block  $N_1$ of the matrix   
  $\times g \in  \C(\widehat{\times z})$, where $g$ is a general linear form,  is of size $2^{n-1}\alpha$ and it consists of 
two fine diagonal blocks of size $2^{n-2}\alpha$ each.    
Note that 
$$
\left\{ 
\begin{array}{l}
U_1:=A/(0:z^2)+(z) \cong K[x_0, x_1, \cdots, x_{n-2}]/(x_0^{\alpha}, x_1 ^2, \cdots, x_{n-2}^2)  \\   
U_2:=(0:z)+(z)/(z) \cong K[x_0, x_1, \cdots, x_{n-2}]/(x_0^{\alpha}, x_1 ^2, \cdots, x_{n-2}^2)[-1]
\end{array}
\right.
$$
Since  $U_1$ and $U_2$ have the SLP,    if $g \in A$ is a general linear form,  then 
the rank of $\times \overline{g^j} \in \End(U_i)$  can be computed from the Hilbert series of $U_i$.  
Now let $G_i$ be a matrix for $\overline{g} \in \End(U_i)$.  
Then  we have 
$\rank (G_1)^3= \dim U_1/(0:g^3)$ and $\rank G_2 = \dim U_2/0:g$.     
Thus, using Lemma~\ref{generalized_binomial} below,  
$$\rank (G_1) ^3 + \rank G_2 + \rank (\times z) = $$
$$(\dim U_1 - s(n-2)-s'(n-2)-s''(n-2)) + (\dim U_2 - s(n-2)) + 2^{n-1}\alpha $$
$$=2^{n}\alpha - s(n) = \mbox{CoSperner number of } A$$
This completes the proof.  (The algebra $A$ in fact has the SLP. This can be proved using 
\cite[Theorem~5.2]{tHjW07}.)  

Verification of the following lemma is left to the reader. 

\begin{lemma}  \label{generalized_binomial}  
Fix a positive integer $\alpha$. Let $A$ be as above.   Define the polynomial $h_n(q)$ by 
\[
h_n(q)=(q^{\alpha -1}+q^{\alpha -2}+ \cdots + q +1)(q+1)^n
\]
Define the integers $s(n), s'(n), s''(n)$  to be the first three terms of the coefficients of 
the polynomial $h_n(q)$  put  in the decreasing order.  
Then we have 
\begin{enumerate}
\item
$h_n(q)$  is the Hilbert series of $A$.  
\item
$h_{n-2}(q)$  is the Hilbert series of $U_1$ and $U_2[1]$.    
\item
$s(n)$ is the Sperner number of $A$ and $s(n-1)$ is the Sperner number of $U_1$ and $U_2$. 
\item
$2^n\alpha -s(n)$ is the CoSperner number of $A$ and $2^{n-1}\alpha -s(n-1)$ is the 
CoSperner number of $U_1$ and $U_2$.   
\item
$s(n)=s(n-1)+s'(n-1)$, for $n \geq 1$. 
\item
$s(n)=2s(n-2)+s'(n-2)+s''(n-2)$, for $n \geq 2$.  
\end{enumerate}
\end{lemma}

Now we prove Example~\ref{yokohama}.
As in the previous example it suffices to prove the WLP for      
$$A:=K[x_0, x_1, \cdots, x_n]/(x_0 ^{\alpha}, x_1^2, x_1x_2, x_2 ^3, x_2x_3, \cdots, 
x_{n-1}x_{n}, x_n^3 )$$  for any positive integer $\alpha$. 

Put $A^{(n)}=A$, $B^{(n-1)}=A^{(n-2)}[z]/(z^2)$.  Then we have the exact sequence:
\[
0 \ra (z) \ra A^{(n)} \ra A^{(n-1)} \ra 0.  
\]
But 
\[
(z)[1] \cong A/(0:z) \cong K[x_0, x_1, \cdots, x_n]/
(x_0 ^{\alpha}, x_1^2, x_1x_2, x_2^3, x_2x_3,  \cdots, x_{n-1}^3, x_{n-1}, x_n^2) \cong  B^{(n-1)}.  
\]
Note that $A^{(n)}$ and $B^{(n)}$ have Hilbert series $$(1+q+ \cdots + a^{\alpha -1})(1+q)^n.$$   
We are going to induct on $n$ so  we assume the SLP for $A^{(n-2)}$ and $B^{(n-2)}$.  
Put $z=x_n$.  Then one sees easily that $\times z$ has the same Jordan decomposition as 
the one treated in Example~\ref{simplified_yokohama}.  
Thus the same proof as Example~\ref{simplified_yokohama} works verbatim in this case.




\end{document}